\documentclass[final, 3p, times]{elsarticle}

\usepackage{amssymb}

\usepackage{amsmath,amsfonts}
\usepackage{algorithmic}
\usepackage{algorithm}
\usepackage{array}
\usepackage{subfig}

\usepackage{graphicx}

\usepackage{bm}
\usepackage{color}
\usepackage{dsfont}
\usepackage{wasysym}
\usepackage{array}
\usepackage{multirow}
\usepackage{optidef}

\usepackage{amsthm}
\theoremstyle{plain}
\newtheorem{theorem}{Theorem}
\newtheorem{lemma}{Lemma}

\newtheorem{corollary}{Corollary}

\DeclareMathOperator{\rank}{rank}

\DeclareMathOperator{\Ima}{Im}
\DeclareMathOperator{\Ree}{Re}

\newcommand{\RR}{\mathbb{R}}
\newcommand{\CC}{\mathbb{C}}
\newcommand{\W}{\mathbf{W}}
\newcommand{\M}{\mathbf{M}}
\renewcommand{\v}{\mathbf{v}}
\newcommand{\Ws}{\mathbf{W}^\star}
\definecolor{Blue}{rgb}{0,0,0}

\journal{Electric Power Systems Research}

\begin{document}

\begin{frontmatter}



\title{Ex Post Conditions for the Exactness of Optimal Power Flow Conic Relaxations}


\author{Jean-Luc Lupien and Antoine Lesage-Landry}

\affiliation{organization={Polytechnique Montreal, Departement of Electrical Engineering},
            addressline={}, 
            city={Montreal},
            postcode={H3T 1J4}, 
            state={QC},
            country={Canada}}

\begin{abstract}
\textcolor{Blue}{
Convex relaxations of the optimal power flow (OPF) problem provide an efficient alternative to solving the intractable alternating current (AC) optimal power flow. The conic subset of OPF convex relaxations, in particular, greatly accelerate resolution while leading to high-quality approximations that are exact in several scenarios. However, the sufficient conditions guaranteeing exactness are stringent, e.g., requiring radial topologies. In this short communication, we present two equivalent \emph{ex post} conditions for the exactness of any conic relaxation of the OPF. These rely on obtaining either a rank-1 voltage matrix or self-coherent cycles.  Instead of relying on sufficient conditions \emph{a~priori}, satisfying one of the presented \emph{ex post} conditions acts as an exactness certificate for the computed solution. The operator can therefore obtain an optimality guarantee when solving a conic relaxation even when \emph{a~priori} exactness requirements are not met. Finally, we present numerical examples from the MATPOWER library where the \emph{ex post} conditions hold even though the exactness sufficient conditions do not, thereby illustrating the use of the conditions.}

\end{abstract}




\end{frontmatter}


\section{Introduction}

\textcolor{Blue}{The optimal power flow (OPF) problem consists of minimizing generation costs in an electric network while ensuring that all loads are adequately supplied and that the transmission physical constraints are met. The alternating current (AC) optimal power flow (\texttt{AC-OPF}), first presented in \cite{carpentier}, is of primary importance to the grid operator for the safe, reliable operation of the system while minimizing operational costs. However, the \texttt{AC-OPF} is, in general, NP-hard \cite{nphard}. For this reason, convex relaxations of the \texttt{AC-OPF} are ubiquitous in the power systems literature \cite{taylor2015convex, coffrin2018convex}.} The widespread adoption of these methods is due to their computational efficiency in providing lower bounds to the \texttt{AC-OPF}. The computational efficiency of convex relaxations also enables the integration of mixed-integer constraints, necessary for unit commitment or network reconfiguration problems~\cite{low2014convex, Survey2019, taylor2015convex}. An important advantage of a convex relaxation over a convex approximation, e.g., the linearized OPF~\cite{taylor2015convex}, is that a relaxation can provide an exact result while taking much less computational resources than the original problem.

A persistent obstacle to the wider adoption of conic relaxations to solving the OPF is that sufficient conditions for optimality are rarely met. \textcolor{Blue}{This means that, for most cases, there is no guarantee that a solution to a convex relaxation of the OPF will be feasible with respect to the original problem \cite{OPFSurvey}.} Instead of relying on these \emph{a priori} conditions, we propose solving a conic relaxation of the OPF problem first, then verifying optimality using an \emph{ex post} condition. \textcolor{Blue}{The possibility of readily obtaining an optimality certificate could potentially justify solving a convex relaxation of the OPF before attempting the nonconvex \texttt{AC-OPF} directly.}

\paragraph{Contribution}
In this short communication, we present two equivalent \emph{ex post} conditions that guarantee that a solution obtained via any conic relaxation of the \texttt{AC-OPF} resulting from a lifted, linear power flow constraint coincides with the solution to the original problem. Specifically, we show that obtaining an optimal voltage matrix that is either rank-$1$ or has self-coherent cycles from a conic relaxation implies that this solution is exact for the original \texttt{AC-OPF} problem. In particular, this result holds for the tight-and-cheap relaxation and the strong, tight-and-cheap relaxation presented in \cite{bingane2018tight}. \textcolor{Blue}{If either condition is respected, the obtained solution will satisfy all \texttt{AC-OPF} constraints and will minimize its objective function. The operator can then directly implement this solution in the grid without having to check for feasibility or running the NP-hard \texttt{AC-OPF}.} We then illustrate the application of these conditions on MATPOWER cases where sufficient exactness conditions are not met.

\section{Conic relaxations of the optimal power flow}
\label{sec:conic_rel}

We consider an electric power system defined as the graph~$\left(\mathcal{N}, \mathcal{L} \right)$ where $\mathcal{N} \subset \mathbb{N}$ is the set of vertices, i.e., buses, and $\mathcal{L} \subset \mathbb{N}\times\mathbb{N}$ is the set of edges, i.e., powerlines. Let $\mathcal{C}$ denote the set of cycles within the graph $\left(\mathcal{N}, \mathcal{L}\right)$. Each element $\mathbf{c}\subset (\mathcal{N}, \mathcal{L})$ of $ \mathcal{C}$ is defined as a non-empty sequence of edges in which only the first and the last vertices are equal.
Let the admittance of line $ij \in \mathcal{L}$ be denoted by $y_{ij} \in \CC$. 

Let $p_i \in \RR$ and $q_i \in \RR$ be, respectively, the active and reactive power injection/absorption at bus $i\in\mathcal{N}$. Similarly, let $p_{ij} \in \RR$ and $q_{ij} \in \RR$ be, respectively, the active and reactive power flow in line $ij \in \mathcal{L}$. Lastly, let $v_i \in \CC$ be the complex voltage phasor at bus $i$. The \texttt{AC-OPF} problem is expressed as:
\begin{mini!}|s|
{}{\sum_{i \in \mathcal{G}} f_i(p_i) \label{eq:opf_obj}}
{\label{eq:opf}}{}
\addConstraint{p_{ij} + \jmath q_{ij} = v_i\left(v_i^* - v_j^*\right)y_{ij}^*,\label{eq:opf_pf}}{\quad}{ij \in \mathcal{L}}
\addConstraint{\underline{v} \leq v_{i} \leq \overline{v},\label{eq:opf_volt}}{\quad}{i \in \mathcal{N}}
\addConstraint{p_i = \sum_{ij \in \mathcal{L}} p_{ij}, {\quad}}{{q_i = \sum_{ij \in \mathcal{L}} q_{ij}},\label{eq:opf_nodalq}}{\quad}{i \in \mathcal{N}}
\addConstraint{\underline{p}_i \leq p_i \leq \overline{p}_i,{\quad}} {\underline{q}_i \leq q_i \leq \overline{q}_i,\label{eq:opf_q}\quad i \in \mathcal{N}}
\addConstraint{p_{ij}^2 + q_{ij}^2 \leq \overline{s}_{ij}^2,\label{eq:opf_apppower}}{\quad}{ij \in \mathcal{L},}
\end{mini!}
where~\eqref{eq:opf_obj} is a convex cost function commonly set as a convex quadratic function~\cite{taylor2015convex},~\eqref{eq:opf_pf} is the power flow constraint,~\eqref{eq:opf_volt} imposes the minimum/maximum ($\underline{v}$ and $\overline{v}$, respectively) voltage magnitude limits,~\eqref{eq:opf_nodalq} are the active and reactive nodal power balance,~\eqref{eq:opf_q} impose the minimum $\underline{p}$ ($\underline{q}$) and maximum $\overline{p}$ ($\overline{q}$) active (reactive) power limits, and~\eqref{eq:opf_apppower} imposes the apparent power limit $\overline{s}_{ij}$ for each line.

A common set of OPF convex relaxations lifts the dimensionality of the voltage variable to rewrite~\eqref{eq:opf_pf} as a linear constraint: 
\begin{equation}
p_{ij} + \jmath q_{ij} = \left(\W_{ii} - \W_{ij}\right)y_{ij}^* \quad ij \in \mathcal{L}, \label{eq:pf_W}
\end{equation}
where $\W \in \CC^{n \times n}$ is a Hermitian matrix, in addition to converting~\eqref{eq:opf_volt} to
\begin{equation}
\underline{v}^2 \leq \W_{ii} \leq \overline{v}^2. \label{eq:opf_new_v}
\end{equation}
If the constraints
\begin{align}
\W &\succeq 0 \label{eq:psd}\\
\rank \W &= 1 \label{eq:rank},
\end{align}
are paired with~\eqref{eq:pf_W} and~\eqref{eq:opf_new_v}, then the resulting OPF problem \textcolor{Blue}{(\texttt{AC-OPF-Eq}), i.e., \eqref{eq:opf_obj} s.t.~\eqref{eq:opf_nodalq}$-$\eqref{eq:opf_apppower},~\eqref{eq:pf_W}$-$\eqref{eq:rank},} is equivalent to~\eqref{eq:opf}~\cite{taylor2015convex}. Common convex relaxations leveraging~\eqref{eq:pf_W} include the semi-definite relaxation (\texttt{OPF-SDR}) and the similar chordal relaxation~\cite{low2014convex}, the second-order cone relaxation (\texttt{OPF-SOCR}), the new semi-definite relaxation (\texttt{OPF-nSDR}), the tight-and-cheap relaxation (\texttt{OPF-TCR}), and the strong tight-and-cheap relaxation (\texttt{OPF-STCR}). 
In \texttt{OPF-SDR}, the non-convex constraint~\eqref{eq:rank} is omitted. Exactness is achieved by \texttt{OPF-SDR} in radial topologies~\cite{taylor2015convex}. The problem \texttt{OPF-SDR} is further relaxed as \texttt{OPF-SOCR} and a non-negativity constraint is applied to the $1\times1$ and $2\times2$ principal minors of $\W$ instead of~\eqref{eq:psd}. In this case, sufficient conditions for exactness exists, viz., radial topologies, infinite line limits, and capability of loads to be over-supplied~\cite[Theorem 3.1]{taylor2015convex}. The \texttt{OPF-nSDR} is similar to \texttt{OPF-SDR} but substitutes~\eqref{eq:psd} and~\eqref{eq:rank} by
\begin{equation}
\W \succeq \mathbf{v} \mathbf{v}^\text{H}, \label{eq:nSDR}
\end{equation}
where $\mathbf{v} \in \CC$ is the vector of voltage phasors, and adds
$\W_{11} \leq \left(\underline{v}_1 + \overline{v}_1 \right) \Ree\left\{ v_1 \right\} - \underline{v}_1 \overline{v}_1$, and
$\Ima \left\{v_1\right\} = 0$,
which are obtained by omitting non-convex constraints resulting from using the reformulation-linearization technique on $\W_{ij}=v_i v_j^*$ \cite{bingane2018tight}. A process similar to what was used to obtain \texttt{OPF-SOCR} from \texttt{OPF-SDR} is utilized to derive \texttt{OPF-TCR} from \texttt{OPF-nSDR}: \eqref{eq:nSDR} is relaxed such that it is only enforced on the $2\times 2$ principal submatrices of $\W$ and the corresponding voltage magnitude vector entries associated with a line $ij \in \mathcal{L}$. 
Lastly, considering that node $1$ is a slack bus, we can use \texttt{OPF-STCR}, which sits between \texttt{OPF-SDR} and \texttt{OPF-TCR} in terms of tightness. In \texttt{OPF-STCR}, the previous constraint is strengthened by directly embedding within it the relations $v_1 = \sqrt{\W_{11}}$ and $v_i = \frac{\W_{i1}}{v_1}$.
This last relaxation is also shown to be equivalent to \texttt{OPF-SDR} in radial topologies~\cite[Proposition 4]{bingane2018tight}. It then follows that \texttt{OPF-STCR} is exact in such topology. Finally, \texttt{OPF-TCR} and \texttt{OPF-STCR} are shown to be tighter than \texttt{OPF-SOCR}~\cite{bingane2018tight}.

\section{\emph{Ex post} exactness conditions}
In this section, we provide two equivalent conditions that can be tested after solving the relaxed OPF problem to check for the exactness of the conic relaxation based on \eqref{eq:pf_W}. If the condition is met, then one can claim that the relaxed solution is equivalent to the original OPF problem's global optimum. The \emph{ex post} condition applies to all relaxations of the \texttt{OPF-SDR} like the ones presented in Section~\ref{sec:conic_rel}. It also trivially applies to \texttt{OPF-SDR}. In this section, we first establish a positive semi-definiteness condition and then state our main result.

\begin{lemma}
Let $\W \in \CC^{n\times n}$ be a Hermitian matrix with $\W_{11} > 0$. If $\rank \W = 1$, then $\W$ is positive semi-definite.
\label{lem:rank}
\end{lemma}

\begin{proof}
Let $\lambda_1\left(\W\right) \leq \lambda_2\left(\W\right) \leq \ldots \leq \lambda_n\left(\W\right)$ be the $n$ real eigenvalues of $\W$. Because $\rank \W = 1$, we have $\lambda_1\left(\W\right) = \lambda_2\left(\W\right) = \ldots = \lambda_{n-1}\left(\W\right) = 0$ and $\lambda_{n}\left(\W\right) \neq 0$.

The Cauchy Interlace Theorem for Hermitian matrix~\cite{cauchyShort, cauchy2004} states that for matrix $\W$ and its $m\times m $ principal submatrix~$\M$ where $m < n$, we have
$\lambda_k\left(\W\right) \leq \lambda_k\left(\M\right) \leq \lambda_{k+n-m}\left(\W\right)$,
for $k=1,2\ldots,m$. Consider the $1\times1$ principal submatrix~$\M = \W_{11}$. Then, the Cauchy Interlace Theorem yields
$\lambda_1\left(\W\right) \leq \W_{11} \leq \lambda_{n}\left(\W\right).$
By assumption, $\W_{11} > 0$ and thus $\lambda_{n}\left(\W\right) > 0$. Because all eigenvalues of $\W$ are non-negative, it follows that $\W \succeq 0$.
\end{proof}
We remark that the assumption of Lemma~\ref{lem:rank} are always met by definition in the OPF problem. Using this lemma, we can provide an \emph{ex post} condition for exactness of a large family of OPF relaxations.


\begin{theorem}\label{thm:exact}
Let \textcolor{Blue}{$\mathcal{W} \subset \mathbb{C}^{n\times n}$ be the set of} optima of a convex relaxation of the optimal power flow problem based on the constraint with lifted dimensionality~\eqref{eq:pf_W}
for all $ij \in \mathcal{L}$. \textcolor{Blue}{The relaxation is exact if and only if there exists $\mathbf{W}^\star \in \mathcal{W}$ such that $\rank \Ws = 1$.} \label{thm:exact_iff}
\end{theorem} 



\begin{proof}
\textcolor{Blue}{We first show the direct implication via its contrapositive. If no rank-1 matrix $\Ws$ exists within $\mathcal{W}$, then no $\mathbf{W} \in \mathcal{W}$ can be feasible with respect to constraint~\eqref{eq:rank} of \texttt{AC-OPF-Eq}, which is equivalent to the original problem~\cite{taylor2015convex}. Therefore, the relaxation cannot be exact and we have established the direct implication.
}

\textcolor{Blue}{Next, we show the converse, i.e., that the existence of a rank-1 matrix $\Ws \in \mathcal{W}$ ensures the relaxation exactness.}
We observe that if $\Ws \succeq 0$ and $\rank \Ws = 1$, then a relaxed OPF using~\eqref{eq:pf_W} is exact because this implies that its solution is such that $\Ws = \v \v^\text{H}$ and, therefore, is feasible with respect to the exact OPF problem.

The voltage bound constraints $\underline{v}^2 \leq \Ws_{ii} \leq \overline{v}^2$ imposes $\Ws_{11} > 0$. By definition, $\Ws$ is Hermitian. Hence, using Lemma~\ref{lem:rank} we have $\Ws \succeq 0$. It follows that a convex relaxation of the OPF based on~\eqref{eq:pf_W} is exact and
$\Ws = \left(\sqrt{\lambda}\mathbf{u}\right)\left(\sqrt{\lambda}\mathbf{u}\right)^\text{H} = \v \v^\text{H}$,
where $\mathbf{u} \in \CC^n$ is the eigenvector of $\Ws$ corresponding to the non-zero eigenvalue $\lambda$. \textcolor{Blue}{Combining the direct implication and the converse completes the proof.}
\end{proof}

Consequently, if the voltage matrix $\Ws$ obtained by solving \texttt{OPF-SOCR}, \texttt{OPF-TCR}, or \texttt{OPF-STCR} to optimality is rank-$1$, then these relaxations are exact. Theorem~\ref{thm:exact} also implies that the simplest and loosest relaxation of that type for which~\eqref{eq:opf_pf} and~\eqref{eq:opf_volt} are substituted by~\eqref{eq:pf_W} and~\eqref{eq:opf_new_v} in~\eqref{eq:opf} is exact if $\rank \W^\star = 1$.

In practice, for most relaxations, the matrix $\Ws$ is sparse because only elements representing power lines are constrained and grids are not fully connected. To verify the rank condition, one must instead check if there exists a completion~$\W^{\text{c}}$ of the matrix $\Ws$ such that $\rank \W^{\text{c}} = 1$. In this case, the following corollary applies.

\begin{corollary} \label{cor:cycles}
\textcolor{Blue}{Let $\mathcal{W} \subset \mathbb{C}^{n\times n}$ be the set of optima of an optimal power flow relaxation based on~\eqref{eq:pf_W}. The relaxation is exact if and only if there exists $\Ws\in\mathcal{W}$ such that the following two conditions are met:} (i)
for every power line $ij\in\mathcal{L}$,
    $\Ws_{ii}\Ws_{jj} = \left|\Ws_{ij}\right|^2$ and (ii) for every cycle $\mathbf{c}\in\mathcal{C}$ for the graph ($\mathcal{N}$, $\mathcal{L}$), $\Ima\left\{\sum\limits_{ij\in \mathbf{c}} \Ws_{ij}\right\} = 0 \mod 2\pi$,
    where $\Ima\{ \cdot \}$ denotes the imaginary part.
\label{cor:comp
}
\end{corollary}

\begin{proof}
\textcolor{Blue}{ The contrapositive of the implication is shown from the feasibility analysis presented in~\cite{low2014convex}, where coherent cycles are necessary for the computed voltages to be physical. Hence, similarly to Theorem~\ref{thm:exact_iff}'s proof, the implication follows from the establishing its contrapositive. The converse is tackled next.} From~\cite[Theorem 3]{low2012}, we have that a sparse matrix $\W^{\text{s}} \textcolor{Blue}{\in \mathcal{W}}$ admits a unique rank-$1$ completion of $\Ws$ if (i) and (ii) hold\textcolor{Blue}{, which is then itself within $\mathcal{W}$.} 
It then follows from Theorem~\ref{thm:exact} that the matrix $\Ws$ is an exact solution to the original OPF problem~\eqref{eq:opf}. \textcolor{Blue}{Showing the implication and the converse statements completes the proof.}
\end{proof}

\textcolor{Blue}{If these \emph{ex post} conditions are not met, the obtained solution will not be feasible with respect to the \texttt{AC-OPF} constraints. This is because obtaining a solution $\W^*$ that does not have self-coherent cycles and is not rank-1 will not produce a set of physically feasible voltages. This follows from Corollary 1 and~\cite{low2014convex}. 
Therefore, to obtain a feasible solution, one must solve the non-convex \texttt{AC-OPF}~\cite{taylor2015convex}. To promote a good quality solution, the infeasible solution obtained from the relaxation can be used as the initial point for the \texttt{AC-OPF}.}

\section{Numerical Example}

To illustrate cases where the \emph{ex post} conditions described in Theorem~\ref{thm:exact}  and Corollary~\ref{cor:cycles} hold, several MATPOWER test cases were solved using conic relaxations. The MATPOWER format and test cases are described in~\cite{Matpower}.  Every test case presented was run using  the $\texttt{OPF-SDR}$, $\texttt{OPF-SOCR}$, $\texttt{OPF-TCR}$, and $\texttt{OPF-STCR}$ and compared to the optimal solution given by MATPOWER's integrated solver. The resolution of the relaxed OPF was achieved using the MOSEK solver \cite{mosek}. The tolerance utilized to determine the rank of a matrix is $\epsilon_1=10^{-6}$ which was used to determine the exactness of the \texttt{OPF-SDR} relaxation. The tolerances used to determine satisfaction of Corollary~\ref{cor:cycles}'s conditions (i) and (ii) are $\epsilon_2=10^{-4}$ and $\epsilon_3=10^{-3}$, respectively. Obtaining a rank-$1$ matrix under this tolerance was used to determine exactness of the \texttt{OPF-SDR} relaxation. Satisfaction of Corollary~\ref{cor:cycles}'s conditions (i) and (ii) is used to determine of the exactness of the \texttt{OPF-SOCR}, \texttt{OPF-TCR}, and \texttt{OPF-STCR} relaxations. The numerical results are presented in Table~\ref{tab:both_tables}.

\begin{table}[htbp]
    \centering
\subfloat[\emph{Ex post} condition satisfaction]{%
    \label{tab:Exact}%
    \begin{tabular}{ccccc}
\hline

\hline
\textbf{Case} \texttt{OPF-} & \texttt{\textbf{SDR}} & \texttt{\textbf{SOCR}} & \texttt{\textbf{TCR}} & \texttt{\textbf{STCR}} \\ \hline
case9 & \checkmark &  & \checkmark &\checkmark \\
case14 & \checkmark &  & \checkmark &\checkmark \\
case22 & \checkmark &  & \checkmark &\checkmark \\
case33bw & \checkmark &  & \checkmark &\checkmark \\
case39 & \checkmark & & \checkmark&\checkmark \\ 
case69 & \checkmark &  & & \checkmark \\
case141 & \checkmark  & & & \checkmark\\
\hline

\hline
\end{tabular}
}
  \qquad\qquad%
\subfloat[Optimality gap]{%
    \label{tab:Looseness}%
\begin{tabular}{ccccc}
\hline

\hline
\textbf{Case} \texttt{OPF-} & \texttt{\textbf{SDR}} & \texttt{\textbf{SOCR}} & \texttt{\textbf{TCR}} & \texttt{\textbf{STCR}} \\ \hline
case9 & 0.00 & 0.00  & 0.00 &0.00 \\
case14 & 0.00 & 0.08  & 0.00 &0.00 \\
case22 & 0.00 & 0.00  & 0.00 & 0.00 \\
case33bw & 0.00 & 0.01 & 0.00 &0.00 \\
case39 & 0.00 & 0.01 & 0.00 & 0.00\\
case69 & 0.00 &0.03 &0.01 & 0.00 \\
case141 & 0.00 & 0.04 &0.01 & 0.00\\
\hline

\hline
\end{tabular}
}
    \caption{Results of Conic Relaxations of the \texttt{AC-OPF}}
    \label{tab:both_tables}
\end{table}

We remark that the optimality gap for the $\texttt{OPF-SOCR}$ is very small for many cases but the requirement of the existence of a rank-$1$ completion is rarely met. This suggests that the $\texttt{OPF-SOCR}$ is a good relaxation if rarely exact. Even if the $\texttt{OPF-TCR}$ takes approximately the same resolution time as $\texttt{OPF-SOCR}$~\cite{bingane2018tight}, it produces many solutions that respect the \emph{ex post} conditions. This is further justification for the adoption of such approaches in future OPF implementations. 

\section{Conclusion}

In this work, two readily-implementable \emph{ex post} conditions guaranteeing that the solution to a conic relaxation of the OPF is exact are established. \textcolor{Blue}{Specifically, obtaining a rank-1 voltage matrix or self-coherent cycles within the voltage matrix indicate that the solution obtained from a conic relaxation coincide with the optimal solution to the \texttt{AC-OPF}.} Numerical simulations of test cases taken from the MATPOWER database are presented to illustrate the application of the conditions. The results notably indicate that the strong, tight-and-cheap relaxation is the most likely relaxation to yield exact solutions after the more costly semi-definite relaxation. \textcolor{Blue}{The good performance of this relaxation as compared to looser relaxations points towards a good compromise between computational complexity and tightness which was also highlighted in \cite{bingane2018tight}. Considering the significant time reduction of running a conic relaxation of the \texttt{OPF} rather than the non-convex \texttt{AC-OPF}, this result suggests that solving the strong, tight-and-cheap relaxation could accelerate the grid operator's resolution speed over long time-periods. The more widespread use of the strong, tight-and-cheap relaxation in future research could increase solution quality over the more prevalent second-order cone relaxation without sacrificing as much computational complexity as a
semi-definite relaxation would. In the future, the scalability of the \emph{ex post} conditions will be considered to reinforce the presented results.}





\bibliographystyle{elsarticle-num}
\bibliography{ref.bib}

\end{document}